\documentclass[journal,transmag]{IEEEtran}
\usepackage{ifpdf}

\usepackage{cite}

\ifCLASSINFOpdf
  \usepackage[pdftex]{graphicx}
\else
   \usepackage[dvips]{graphicx}
\fi

\usepackage{amsmath}
\usepackage{amssymb}
\usepackage{mathtools}

\usepackage{algorithmic}

\usepackage{array}

\ifCLASSOPTIONcompsoc
 \usepackage[caption=false,font=normalsize,labelfont=sf,textfont=sf]{subfig}
\else
 \usepackage[caption=false,font=footnotesize]{subfig}
\fi

\usepackage{stfloats}

\ifCLASSOPTIONcaptionsoff
 \usepackage[nomarkers]{endfloat}
\let\MYoriglatexcaption\caption
\renewcommand{\caption}[2][\relax]{\MYoriglatexcaption[#2]{#2}}
\fi

\usepackage{url}
\usepackage{hyperref}

\usepackage{eso-pic}
\AddToShipoutPicture*{\footnotesize\sffamily\raisebox{0.75cm}{\hspace{1.35cm}\fbox{\parbox{\textwidth}{\copyright\,2023 IEEE.  Personal use of this material is permitted.  Permission from IEEE must be obtained for all other uses, in any current or future media, including reprinting/republishing this material for advertising or promotional purposes, creating new collective works, for resale or redistribution to servers or lists, or reuse of any copyrighted component of this work in other works.}}}}

\begin{document}


\title{2D Eddy Current Boundary Value Problems for \\ Power Cables with Helicoidal Symmetry}


\author{\IEEEauthorblockN{Albert Piwonski\IEEEauthorrefmark{1},
Julien Dular\IEEEauthorrefmark{2},
Rodrigo Silva Rezende\IEEEauthorrefmark{1}, 
Rolf Schuhmann\IEEEauthorrefmark{1}}
\IEEEauthorblockA{\IEEEauthorrefmark{1}Theoretische Elektrotechnik,
Technische Universität Berlin, Berlin, Germany}
\IEEEauthorblockA{\IEEEauthorrefmark{2}Applied and Computational Electromagnetics,  University of Liège, Liege, Belgium}%
}


\markboth{CEFC 2022: P09: Static and Quasi Static Fields No. 2, MAGCON-22-11-0332.R1, Conference Paper Number: 154}%
{Shell \MakeLowercase{\textit{et al.}}: Bare Demo of IEEEtran.cls for IEEE Transactions on Magnetics Journals}


\IEEEtitleabstractindextext{%
\begin{abstract}

Power cables have complex geometries in order to reduce their AC resistance. The cross-section of a cable consists of
several conductors that are electrically insulated from each other to counteract the current displacement caused by 
the skin effect. Furthermore, the individual conductors are twisted over the cable’s length. This geometry has a non-standard symmetry -- a combination of translation and rotation. Exploiting this property allows formulating a dimensionally reduced boundary value problem. Dimension reduction is desirable, otherwise the electromagnetic modeling of these cables becomes impracticable due to tremendous computational efforts. We investigate 2D eddy current boundary value problems which still allow the analysis of 3D effects, such as the twisting of conductor layers.

\end{abstract}

}

\maketitle
\IEEEdisplaynontitleabstractindextext
\IEEEpeerreviewmaketitle


\section{Introduction}
\label{sec:introduction}
\begin{figure}[ht]
    \centering
    \includegraphics[width=0.4\textwidth]{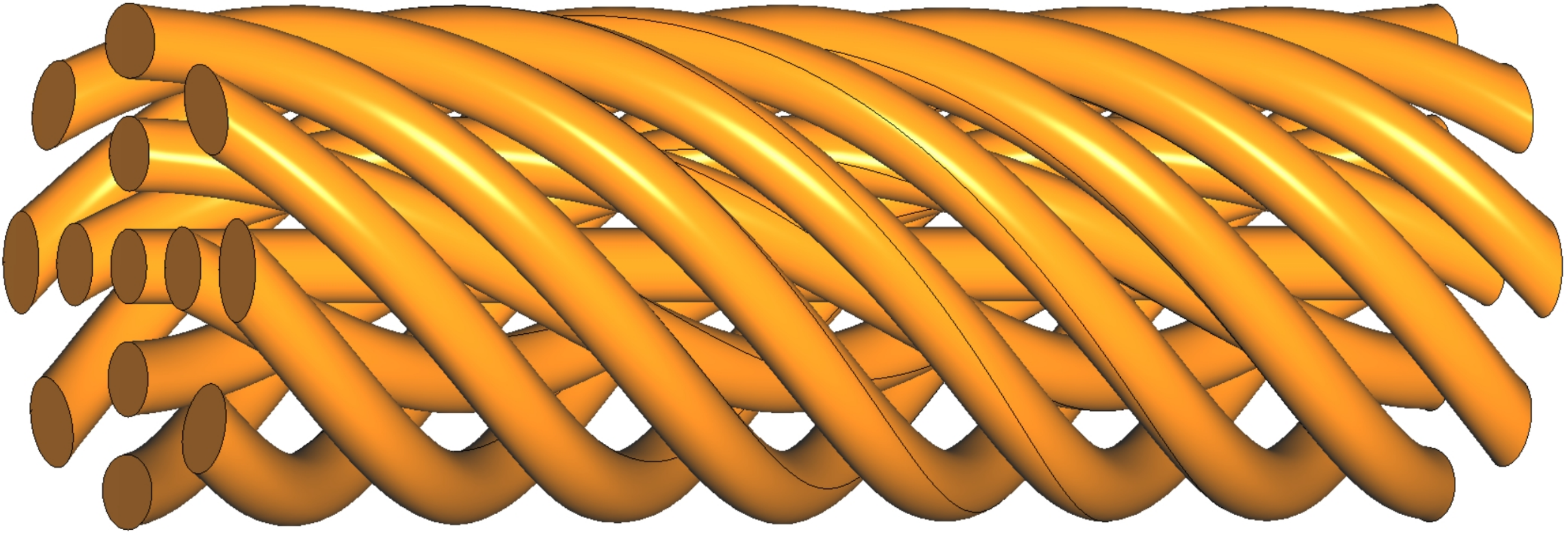}
    \caption{Generic model of a cable's inner conductor: Modeled using \cite{a1}.}
    \label{fig:cst_cable}
\end{figure}

\IEEEPARstart{P}{ower} cables are important elements in the transmission chain of electric power from generator to consumer. Special cable designs are used for AC operation to counteract the undesirable current displacement caused by the skin and proximity effect. Although there are many different cable designs, most have in common that their inner conductors' cross-section is divided into several conductors, which are twisted and electrically insulated from each other (see fig. \ref{fig:cst_cable}).

Analyzing the impedance of power cables is a challenging task, both in laboratories and using numerical simulations. Real-world measurements either are expensive and time-consuming, because a thermal steady state of the cable has to be achieved (calorimetric methods), or they are sensitive to interference signals coupling into the voltage pick-up loop (electrical methods), s.t. complex measurement setups are required \cite{a2}. Likewise, solving an eddy current boundary value problem (BVP) numerically in 3D, which models the cable's electromagnetic behaviour in the magnetoquasistatic limit, leads to tremendous computational efforts, due to the complex geometries. Here, a particular challenge is the multiscale problem that arises from modeling the thin insulations of the individual conductors. 

However, if the BVP has a symmetry, computational costs can be scaled down significantly by means of dimension reduction. Choosing proper boundary conditions, an eddy current BVP posed on a domain as in fig. \ref{fig:cst_cable} has what is called a helicoidal symmetry. This symmetry was exploited before for wave propagation and hysteresis loss problems~\cite{a3},\cite{a4}. In contrast to applying periodic boundary conditions, here the model can be solved in 2D.

Our work focuses on using this approach for the time-harmonic analysis of power cables: In sec. \ref{sec:intuition_for_helicoidal_symmetry} we give an intuition for helicoidal symmetries and introduce the coordinate system in which we perform the dimension reduction. Sec.~\ref{sec:2D_eddy_current_BVP_in_helicoidal_coordinates} is dedicated to the finite element formulation and implementation details of the 2D model. In sec. \ref{sec:results} we compare the results with a 3D reference model. 


\section{Intuition for helicoidal symmetry}
\label{sec:intuition_for_helicoidal_symmetry}

Loosely speaking, symmetry means the property of an object to remain the same under geometric transformations. By this is meant here that every cross-section perpendicular to the longitudinal direction of the cable looks the same after the function composition of rotation and translation. In order to exploit symmetries for the dimension reduction of BVPs, it is desirable to use a coordinate system in which the object being studied appears to be the same in one direction. For objects that appear helical in Cartesian coordinates the helicoidal coordinate system fulfills that requirement. 

\subsection{Helicoidal coordinates}
\label{subsec:helicoidal_coordinates}

\begin{figure}[h!]
    \centering
    \includegraphics[width=0.48\textwidth]{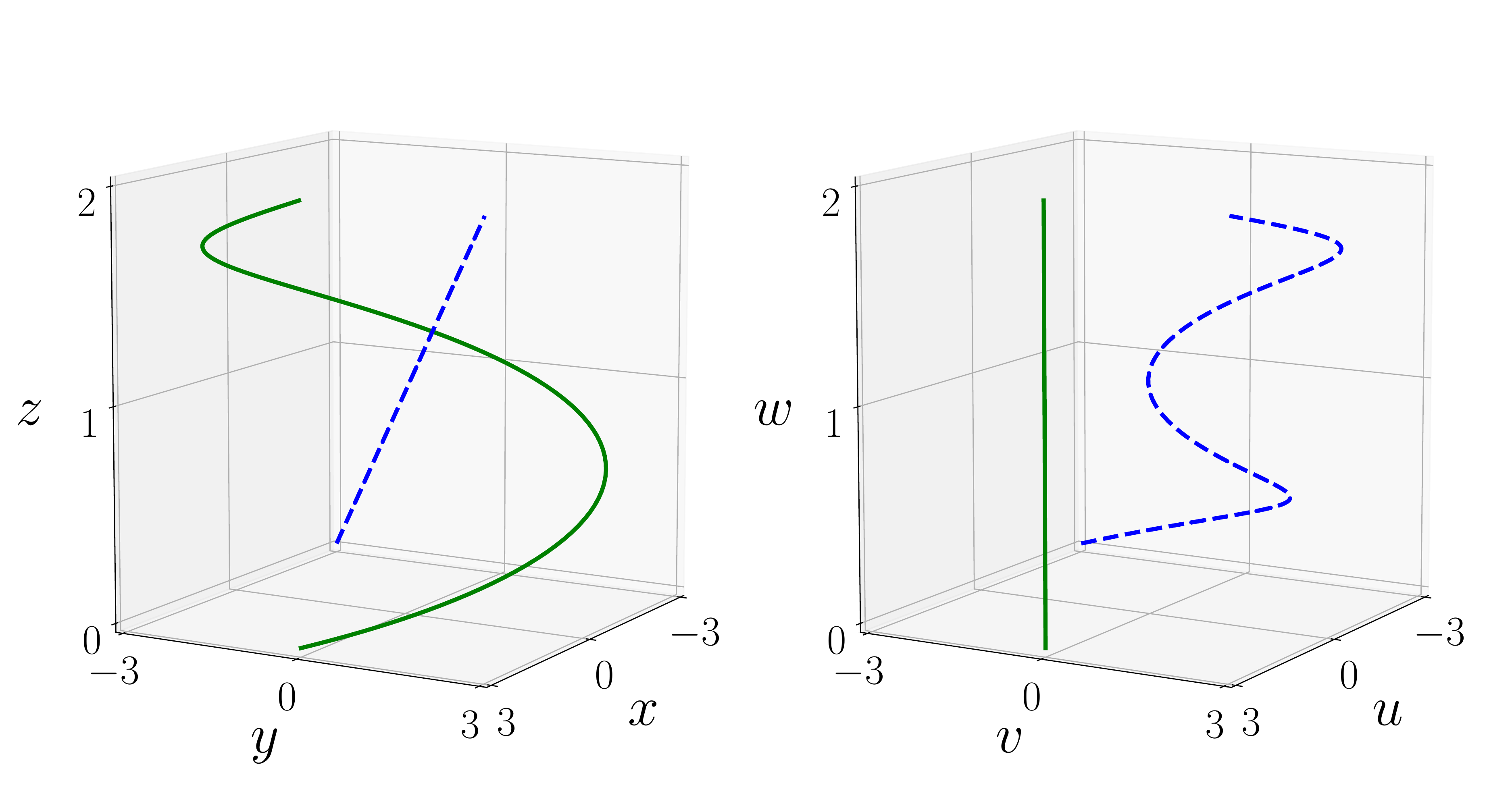}
    \caption{Geometric objects represented in Cartesian (left, arbitrary units) \& helicoidal coordinates (right, same arbitrary units): Helical objects appear straight, which is not necessarily true for straight lines (e.g., see the dashed line between $[-3,\,-3,\,0]^\top$ and $[3,\,3,\,2]^\top$). Further, note that the $(u,\,v,\,w)$ coordinate system must be understood here as non-directional since it is non-orthogonal ($w$-axis is not orthogonal to the $uv$-plane).}
    \label{fig:understanding_transforamtion}
\end{figure}

In the following, we denote points represented in the Cartesian coordinate system $(x,\,y,\,z)$ as $\mathbf{p}_{xyz}\coloneqq [x,\,y,\,z]^\top$, whereas points represented in the helicoidal coordinate system $(u,\,v,\,w)$ are denoted as $\mathbf{p}_{uvw}\coloneqq [u,\,v,\,w]^\top$. The change of coordinates is achieved by the map $\boldsymbol\phi: \Omega_{xyz} \rightarrow \Omega_{uvw}$ and its inverse $\boldsymbol\phi^{-1}:\Omega_{uvw} \rightarrow \Omega_{xyz}$, where $\Omega_{xyz},\, \Omega_{uvw} \subset \mathbb{R}^3$: 

\begin{align}
    \boldsymbol\phi(\mathbf{p}_{xyz}) &= \mathbf{p}_{uvw} = 
    \begin{bmatrix}
        + x \cos(z\alpha/\beta) + y \sin(z\alpha/\beta) \\
        - x \sin(z\alpha/\beta) + y \cos(z\alpha/\beta) \\
        +z
    \end{bmatrix}, \label{eq:xyz_2_uvw}
    \end{align}
    \begin{align}
    \boldsymbol\phi^{-1}(\mathbf{p}_{uvw}) &= \mathbf{p}_{xyz} =
    \begin{bmatrix}
        +u\cos(w\alpha/\beta) -v\sin(w\alpha/\beta) \\
        +u\sin(w\alpha/\beta) +v\cos(w\alpha/\beta) \\
        +w
    \end{bmatrix}\label{eq:uvw_2_xyz}.
\end{align}
Here, parameters $\alpha$, $\beta$ are related to the number of turns and to the total longitudinal length of the helical object of interest, i.e., for different geometries $\boldsymbol\phi$ is defined differently as well. The effect of the (global) coordinate transformation is demonstrated in fig. \ref{fig:understanding_transforamtion}.
\vspace{-0.45cm}

\subsection{Symmetry cell derivation via the theory of envelopes}
\label{subsec:symmetry_cell_derivation_via_theory_of_envelopes}

If a BVP has a symmetry on a domain $\Omega$, then it is sufficient to solve a lower-dimensional BVP on a lower-dimensional subset of $\Omega$, called the symmetry cell. Then the solution on the whole domain $\Omega$ can be derived using symmetry transformations \cite{a4}. Assume an eddy current BVP posed on the domain as in fig. \ref{fig:cst_cable} supplemented with an additional concentrically arranged and perfect electrically conductive (PEC) cylinder that models the cable's shielding. Then, a cross-section generated by a cut orthogonal to the longitudinal direction of the cable qualifies for a symmetry cell. Here it should be noted that due to the twist, the conductors' cross-sections are not circular anymore. 

In the following, we derive the symmetry cell using the mathematical theory of envelopes \cite{c2}. Assume the parametrization $\boldsymbol\gamma$ of a helix curve in Cartesian coordinates with a clockwise rotation:
\begin{align}
    \boldsymbol\gamma: [0,\,2\pi] &\rightarrow \mathbb{R}^3, \, t \mapsto [\underbrace{r\cos(\alpha t)}_{\coloneqq\gamma_x(t)},\, \underbrace{r\sin(\alpha t)}_{\coloneqq \gamma_y(t)},\, \underbrace{\beta t}_{\coloneqq\gamma_z(t)}]^\top, \label{eq:helix_curve}
\end{align}
where $t$ is the parametrization parameter, $\alpha\in\mathbb{R}^+$ is the number of turns per $2\pi$, $\beta\in\mathbb{R}^+$ is the longitudinal length divided by $2\pi$ and $r$ is the helix's radius. Then, sweeping a sphere with some positive radius $r_c$ (conductor radius) along~$\boldsymbol\gamma$, allowing different radii and shifted trigonometric functions in~(\ref{eq:helix_curve}), leads to the helicoidal symmetric structure in fig.~\ref{fig:cst_cable}. We are now interested in calculating the exact shape of the cross-section at $z=0$. Although other choices for the fixed $z$-coordinate would also lead to valid symmetry cells, this choice is easier, because there, $x=u, \,y=v, \, z=w=0$ holds (see eq.~(\ref{eq:xyz_2_uvw})). 

\begin{figure}[ht]
    \centering
     \includegraphics[width=0.3\textwidth]{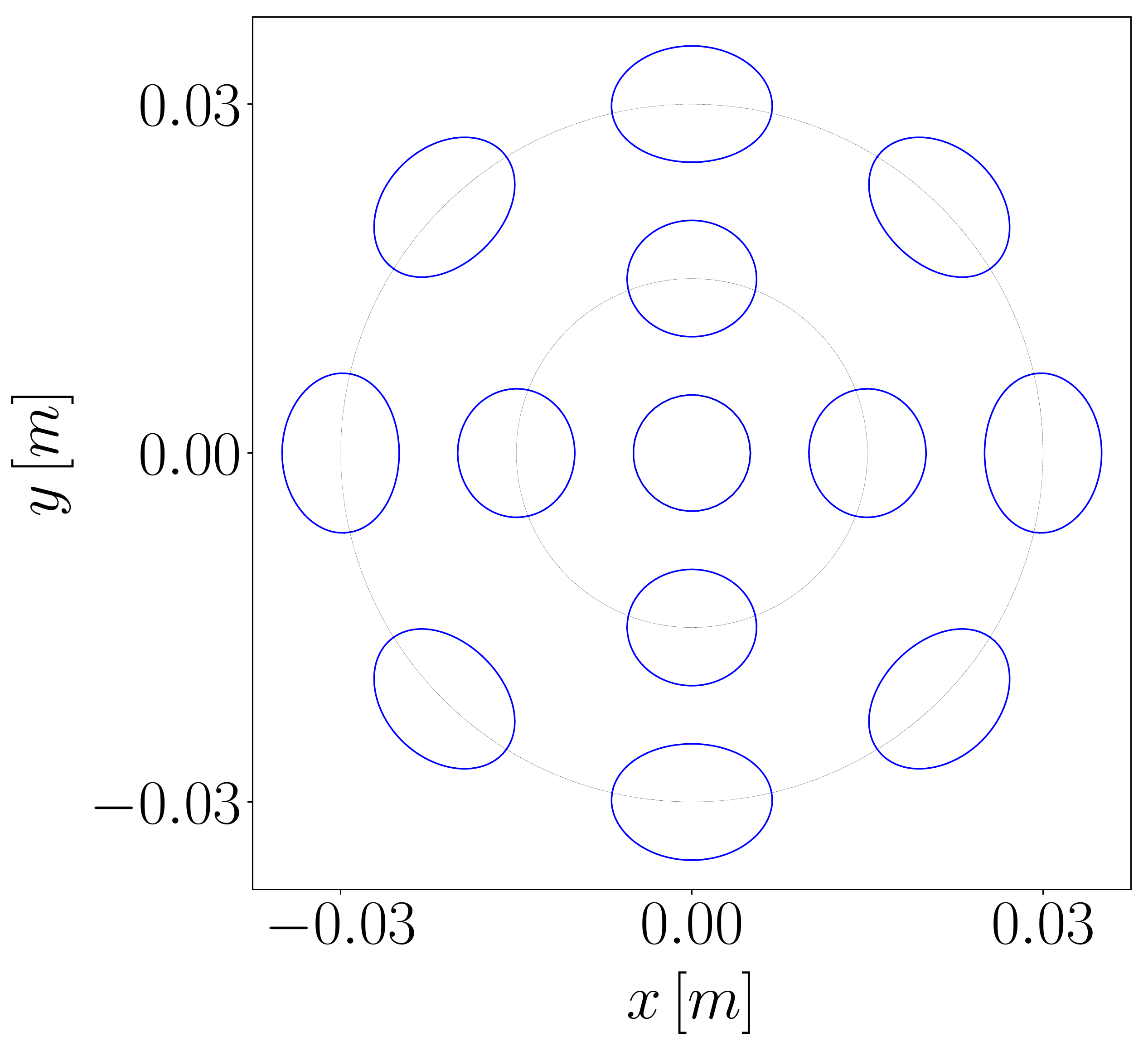}
     \caption{Symmetry cell computed exactly via the theory of envelopes (corresponds to the front view 
     of the cable in fig.~\ref{fig:cst_cable}) for parameter values: $\alpha=1$, $\beta = 0.0318\,\mathrm{m}$, $r_c = 0.5\,\mathrm{cm}$, conductor layers at radii $\{0.0,\,1.5,\,3.0\}\,\mathrm{cm}$.}
     \label{fig:cross_section}
\end{figure}

With this choice, the procedure in order to find the cross-section for a single conductor is as follows: First, a differentiable function $g:[0,\,2\pi]\times\mathbb{R}^2 \rightarrow \mathbb{R}$ is defined as: 
\begin{align}
    g(t,x,y) \coloneqq \left\|[x-\gamma_x(t),\, y-\gamma_y(t), \, 0 - \gamma_z(t)]^\top\right\|^2- r_c^2 \label{eq:points_on_ball},
\end{align}
whose zeros describe implicitly the points lying on the sphere with radius $r_c$ centered at some point on the helix curve, i.e., at some value $t$ in the parametrization~(\ref{eq:helix_curve}).
Solving $g(t,x,y) =0$ for $y(t,x)$ defines a family of curves (one curve $y(t,x)$ for each $t\in[0,\,2\pi]$). The envelope of this family of curves is then formed by points fulfilling additionally:
\begin{align}
    \frac{\partial g}{\partial t}(t,x,y) = 0. \label{eq:tangential_condition}
\end{align}
To find these points we replace $y$ by $y(t,x)$ in eq.~(\ref{eq:tangential_condition}) which leads to an analytically solvable root-finding problem, for each considered $t\in[0,\,2\pi]$. The resulting cross-section for an exemplary helicoidal symmetric structure is shown in fig.~\ref{fig:cross_section}. 

\section{2D eddy current BVP in helicoidal coordinates}
\label{sec:2D_eddy_current_BVP_in_helicoidal_coordinates}

\subsection{Weak formulation}
\label{subsec:weak_formulation}

Our approach is based on the $\mathbf{H}$-$\varphi$-formulation for posing the eddy current BVP in the time-harmonic case, which enables a convenient way to impose global conditions, e.g., fixing the total current $I_i$ flowing through the $i$-th conductor~\cite{d0}.

In the Cartesian domain $\Omega_{xyz}$, with conducting subdomain $\Omega_{xyz,c}$, and in frequency domain with frequency $\omega/2\pi$, the $\mathbf{H}$-$\varphi$-formulation reads~\cite{d0}:
\begin{align}
    &\int_{\Omega_{xyz}} \left(j\omega \mu_{xyz} \mathbf{H}_{xyz}\right) \cdot \mathbf{H}_{xyz}'\,\mathrm{d}V\,\nonumber\\
    &\qquad +\int_{\Omega_{xyz,c}} \left(\rho_{xyz}\, \mathrm{curl}\,\mathbf{H}_{xyz}\right) \cdot \mathrm{curl}\, \mathbf{H}_{xyz}'\,\mathrm{d}V = 0,\label{eq:weak_form_basic} 
\end{align}
for the magnetic field $\mathbf{H}_{xyz}$, magnetic permeability $\mu_{xyz}$, electric resistivity~$\rho_{xyz}$, and with trial and test functions $\mathbf{H}_{xyz}'$ to be chosen in appropriate function spaces. The magnetic field~$\mathbf{H}_{xyz}$ is a one-form and is curl-free in the non-conducting domain $\Omega_{xyz,i}$. So far, the formulation is still fully 3D, as the geometry is not $z$-invariant.

We then change the variables into the helicoidal coordinate system. The one-form $\mathbf{H}_{xyz}$, and two-form $\mathrm{curl}\, \mathbf{H}_{xyz}$ (current density $\mathbf{J}_{xyz}$), transform as follows~\cite{c3,lindell}:
\begin{align}
    \mathbf{H}_{xyz}(\mathbf{p}_{xyz}) &= J_{\boldsymbol\phi^{-1}}^{-\top}\ \mathbf{H}_{uvw}(\boldsymbol\phi(\mathbf{p}_{xyz})),\label{eq:H_uvw_2_H_xyz}\\
    \ \mathrm{curl}\,\mathbf{H}_{xyz}(\mathbf{p}_{xyz}) &= \frac{J_{\boldsymbol\phi^{-1}}}{\mathrm{det}(J_{\boldsymbol\phi^{-1}})}\ \mathrm{curl}\,  \mathbf{H}_{uvw}(\boldsymbol\phi(\mathbf{p}_{xyz})),\label{eq:J_uvw_2_J_xyz}
\end{align}
where $J_{\boldsymbol\phi^{-1}}$ denotes the Jacobian of $\boldsymbol\phi^{-1}$ evaluated at point $\mathbf{p}_{uvw}$. Changing variables also introduces a factor $\mathrm{det}(J_{\boldsymbol\phi^{-1}}) = 1$ in the volume integrals in eq.~\eqref{eq:weak_form_basic}. In terms of the helicoidal coordinates, introducing the trial and test function spaces for the magnetic field $V(\Omega_{uvw})$ and $V_0(\Omega_{uvw})$, respectively, we can therefore rewrite this formulation as follows. Seek $\mathbf{H}_{uvw} \in V(\Omega_{uvw})$, such that $\forall\,\mathbf{H}_{uvw}' \in V_0(\Omega_{uvw})$:
\begin{align}
    &\int_{\Omega_{uvw}} \left(j\omega \boldsymbol\mu_{uvw} \mathbf{H}_{uvw}\right) \cdot \mathbf{H}_{uvw}'\,\mathrm{d}V\,\nonumber\\ 
    &\qquad + \int_{\Omega_{uvw,c}} \left(\boldsymbol\rho_{uvw}\, \mathrm{curl}\,\mathbf{H}_{uvw}\right) \cdot \mathrm{curl}\, \mathbf{H}_{uvw}'\,\mathrm{d}V = 0, \label{eq:weak_form}
\end{align}
where the effect of the change of variables is fully contained in two anisotropic material parameters, written as tensors:
\begin{align}
    \boldsymbol\mu_{uvw}(\mathbf{p}_{uvw}) &= \mu_{xyz}(\boldsymbol\phi^{-1}(\mathbf{p}_{uvw})) J_{\boldsymbol\phi^{-1}}^{-1}J_{\boldsymbol\phi^{-1}}^{-\top}\mathrm{det}(J_{\boldsymbol\phi^{-1}}), \label{mu_xyz_2_uvw} \\
    \boldsymbol\rho_{uvw}(\mathbf{p}_{uvw}) &= \rho_{xyz}(\boldsymbol\phi^{-1}(\mathbf{p}_{uvw})) J_{\boldsymbol\phi^{-1}}^{\top}J_{\boldsymbol\phi^{-1}}/\mathrm{det}(J_{\boldsymbol\phi^{-1}}). \label{rho_xyz_2_uvw} 
\end{align}
The key point of the approach is that the product $J_{\boldsymbol\phi^{-1}}^{\top}J_{\boldsymbol\phi^{-1}}$ (and its inverse) is independent of $w$. Since neither the material nor the symmetry cell's shape depend on~$w$, partial derivatives with respect to $w$ of electromagnetic field quantities vanish ($\partial_w\cdot=0$). This allows the dimension reduction of the BVP to 2D, which we solve on the $uv$-plane at $w=0$.

\subsection{Space discretization}
\label{subsec:space_discretization}

The BVP is 2D, but the magnetic field $\mathbf{H}_{uvw} \in V(\Omega_{uvw})$ still has three components. In practice, we treat the in-plane components ($H_u$, $H_v$) and the $w$-component ($H_w$) separately. In the non-conducting domain $\Omega_{uvw,i}$, we want the magnetic field to be curl-free. From eq.~\eqref{eq:J_uvw_2_J_xyz}, because the change of variables is regular, the curl-free condition in $(x,\,y,\,z)$ coordinates translates into:
\begin{align}
\mathrm{curl}\,\mathbf{H}_{uvw} \stackrel{\partial_w\cdot = 0}{=} \begin{bmatrix}
\partial_{v}H_{w}\\
-\partial_{u}H_{w}\\
\partial_{u}H_{v} - \partial_{v}H_{u}
\end{bmatrix} = \mathbf{0}.\label{eq:curlFreeInHelicoidal}
\end{align}

The vector made up by the in-plane components $H_u$ \& $H_v$ can therefore be discretized using 2D Whitney edge functions in the conducting domain $\Omega_{uvw,c}$, and by curl-free functions in the complementary non-conducting domain $\Omega_{uvw,i}$. The curl-free space is spanned by gradients of a scalar potential $\varphi_{uvw}$, plus a basis of the cohomology space~$\mathcal{H}^1(\Omega_{uvw,i})$~\cite{d0,d3}. This defines the function space $V_{uv}(\Omega_{uvw})$.

The component $H_w$ is expressed via node functions in $\Omega_{uvw,c}$, and, using eq.~\eqref{eq:curlFreeInHelicoidal}, is constant in each connected part of $\Omega_{uvw,i}$. This constant will be chosen to fix the axial field on the PEC boundary. This defines the function space $V_{w}(\Omega_{uvw})$. In total, we have $V(\Omega_{uvw}) = V_{uv}(\Omega_{uvw}) \oplus V_{w}(\Omega_{uvw})$. Further, the test function space $V_0(\Omega_{uvw})$ is the same space but with homogeneous essential boundary conditions.

\subsection{Homology and cohomology basis selection}
\label{subsec:homology_and_cohomology_basis_selection}

We solve the BVP (\ref{eq:weak_form}) by using the open-source finite element software GetDP~\cite{d1}, which allows for flexible function space definitions, whereas the meshing process is performed by Gmsh controlled via the Julia API~\cite{d2}, \cite{c4}. 

Using cohomology basis functions, it is possible to impose global conditions of the BVP directly in $V(\Omega_{uvw})$. Gmsh provides an integrated (co)homology solver which outputs a basis for the first homology space $\mathcal{H}_1(\Omega_{uvw,i})$ and cohomology space $\mathcal{H}^1(\Omega_{uvw,i})$~\cite{d3}. In short, homology describes the tunnels through the non-conducting domain (e.g., generated by piercing conductors), whereas cohomology assigns global quantities to them (total current). In our application it is desirable to fix each conductor's total current separately, i.e., each basis element of $\mathcal{H}_1(\Omega_{uvw,i})$ has to capture exactly one disjoint tunnel. 

\begin{figure}[ht]
    \begin{minipage}[c]{0.24\textwidth}
    \centering
    \includegraphics[width=1.0\textwidth]{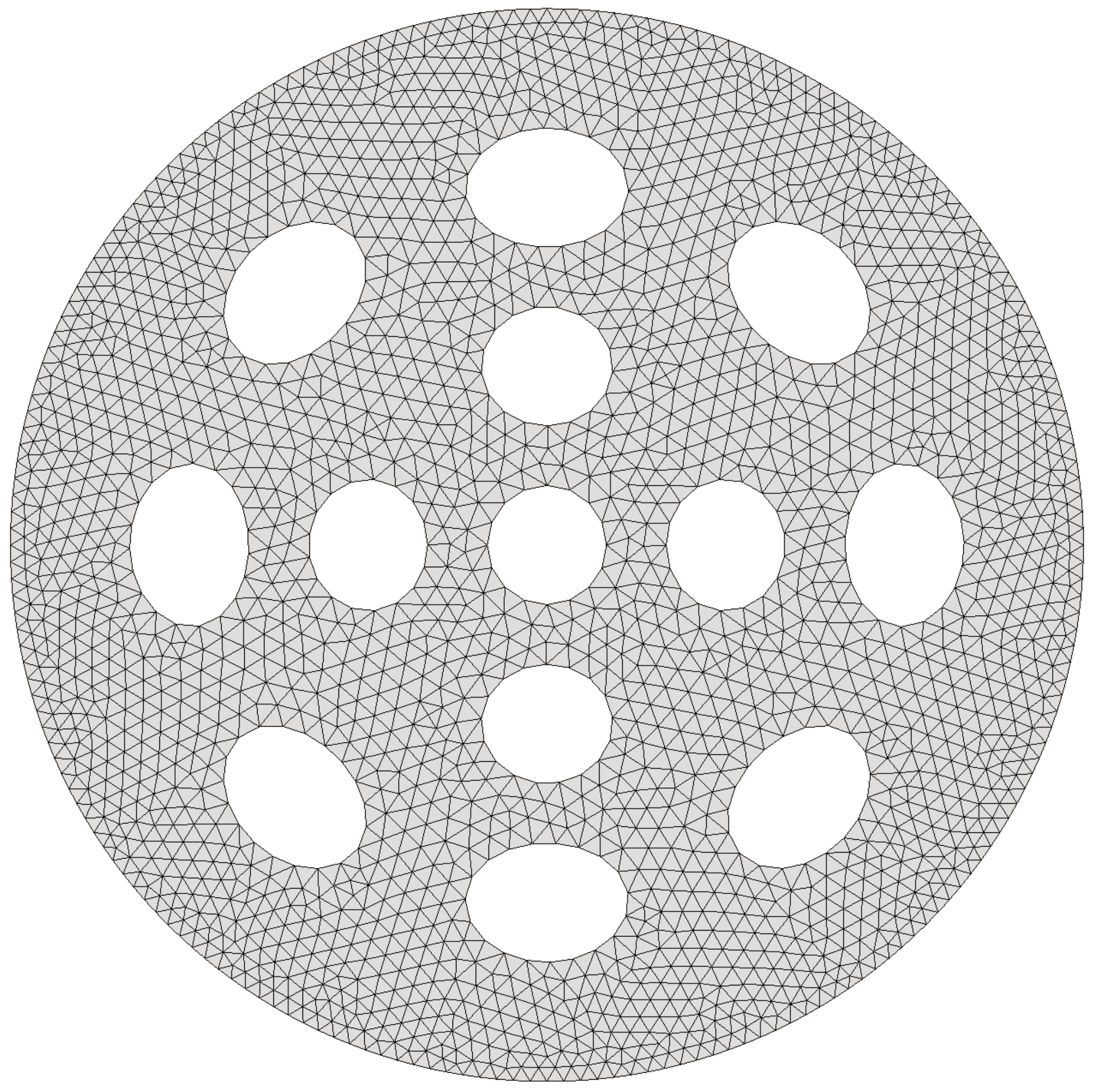}
    \end{minipage}
    \hfill
    \begin{minipage}[c]{0.24\textwidth}
    \centering
    \includegraphics[width=1.0\textwidth]{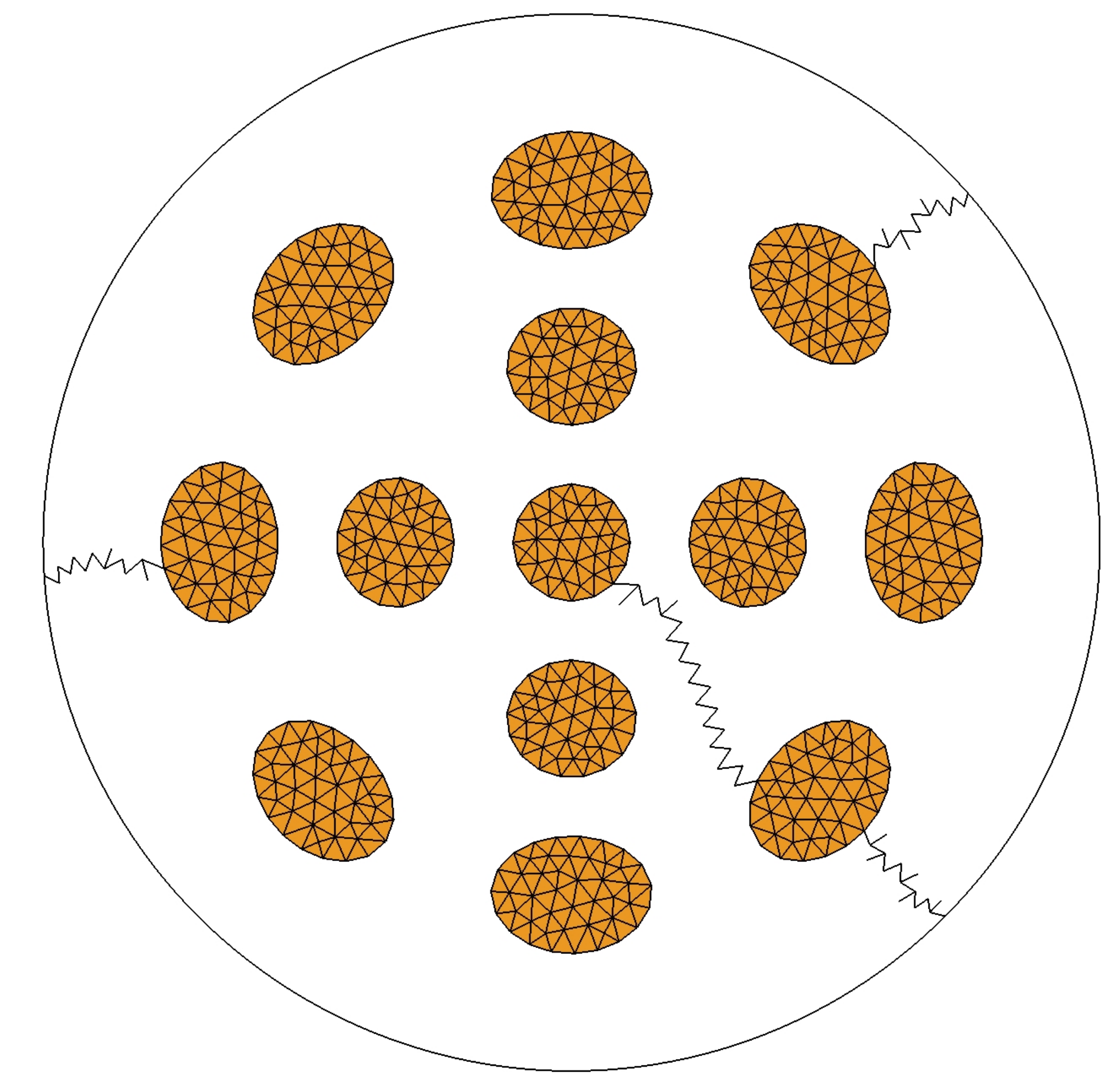}
    \end{minipage}
    \caption{Left: Coarsely triangulated non-conducting domain~$\Omega_{uvw,i}$ with $13$ tunnels/conductors, right: representation of $3$ basis elements spanning parts of $\mathcal{H}^1(\Omega_{uvw,i})$. Due to the basis manipulation, the support of each cohomology basis function (so called thick cut) ranks from a disjoint conductor to $\partial\Omega_{uvw}$.}
    \label{fig:sorted_cohomology}
\end{figure}

This is achieved by first computing a basis for $\mathcal{H}_1(\partial\Omega_{uvw,c})$, where $\partial\Omega_{uvw,c}$ is the boundary of $\Omega_{uvw,c}$. This space is, from a topological perspective, 
indistinguishable from $\mathcal{H}_1(\Omega_{uvw,i})$. Then, a pre-computed, but unwanted basis for $\mathcal{H}^1(\Omega_{uvw,i})$ can be made compatible with the found basis for $\mathcal{H}_1(\Omega_{uvw,i})$ by simple matrix manipulations of the bases representing matrices~(see fig.~\ref{fig:sorted_cohomology})~\cite{d3}.

\vspace{-0.2cm}


\section{Results}
\label{sec:results}

\subsection{3D reference model}
\label{subsec:3D_reference_model}

\begin{figure}[h!]
    \centering
    \includegraphics[width=0.3\textwidth]{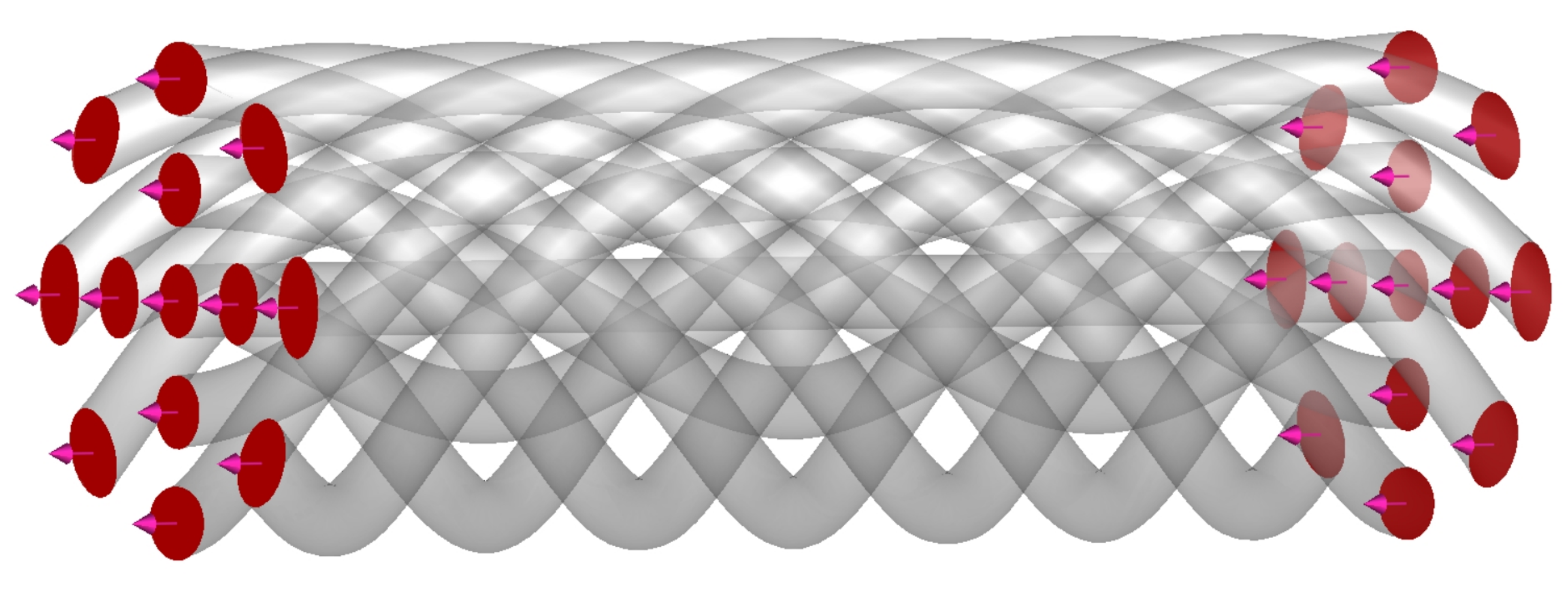}
    \caption{Current excitation ports of the 3D cable model.}
    \label{fig:cst_excitation}
\end{figure}

As a reference, we implemented a 3D cable model with also $13$ helical conductors with a longitudinal length of $0.2\,\mathrm{m}$ in the commercial software CST Studio Suite~\cite{a1}, referred to as CST. In both models, we considered annealed copper for the conductors' material (resistivity $\approx 17.2 \,\mathrm{p}\Omega\mathrm{m}$) and further assumed a non-magnetic material in the whole domain ($\mu_0 = 4\pi\cdot 10^{-7}\, \mathrm{H}/\mathrm{m}$). Further, in both models each conductor carries a total current of amplitude $\sqrt{2}/13 \, \mathrm{A}$ at $f = 50 \, \mathrm{Hz}$. The current constraint is implemented in the 3D model via current ports located at the cuboid bounding box of the computational domain (see ﬁg.~\ref{fig:cst_excitation}). Therefore, the imposed current density bends towards the helical geometries of the conductors only shortly after the ports. To minimize this impact on the results, we considered local field quantities only at the cable's center. In the following, we compare a 2D model discretized by $42.47\,\mathrm{k}$ triangles with a 3D model discretized by $1.19\,\mathrm{M}$ tetrahedra. 
In both models we used first order finite elements resulting in $39.9\,\mathrm{k}$ degrees of freedom in the 2D model and $1.21\,\mathrm{M}$ degrees of freedom in the 3D model.

\subsection{Comparison of local \& global quantities}
\label{subsec:comparison_of_local_and_global_quantites}

The 2D model outputs the finite element approximated magnetic field $\mathbf{H}_{uvw}$ and the current density $\mathbf{J}_{uvw}$, as a secondary quantity, which are then transformed back into Cartesian coordinates using eq.~\eqref{eq:H_uvw_2_H_xyz} and \eqref{eq:J_uvw_2_J_xyz}. 

\begin{figure}[!ht]
    \centering
    \includegraphics[width=0.5\textwidth]{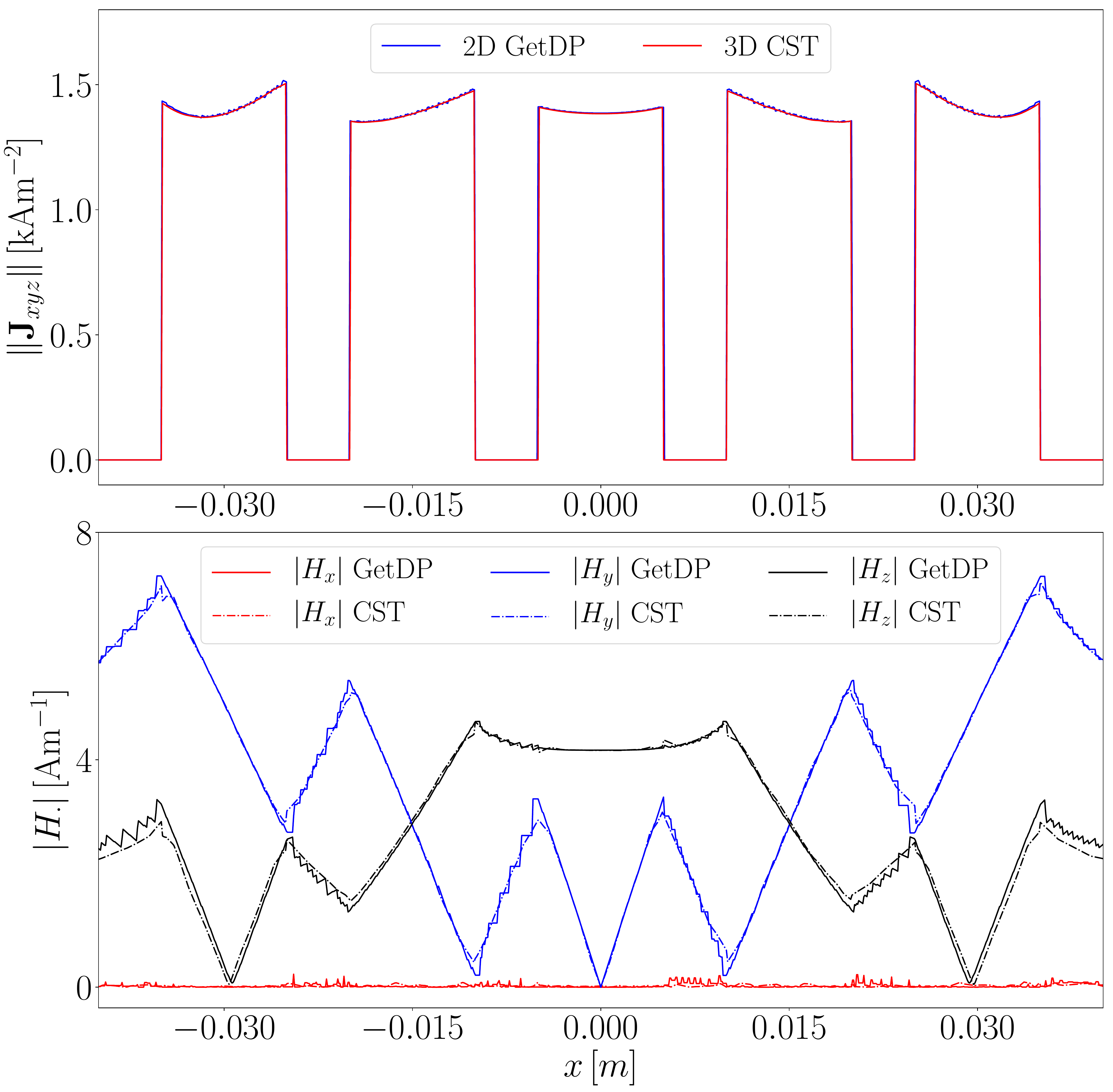}
    \caption{Local comparison: Absolute of current density $\mathbf{J}_{xyz}$ and of magnetic field components $H_x,\, H_y,\,H_z$ along $x$-line. }
    \label{fig:J_Get_DP_vs_CST}
\end{figure}

As a local comparison, $\mathbf{H}_{xyz}$ and $\mathbf{J}_{xyz}$ are evaluated along the $x$-axis. The results depicted in fig.~\ref{fig:J_Get_DP_vs_CST} show a good agreement between both models. The linear finite elements used as ansatz functions for interpolating $\mathbf{H}_{uvw}$ are leading to an element-wise constant current density $\mathbf{J}_{uvw}$, since both quantities are coupled via the curl-operator. However, $\mathbf{J}_{uvw}$ appears jagged represented as $\mathbf{J}_{xyz}$ in the Cartesian coordinate system due to formula~(\ref{eq:J_uvw_2_J_xyz}). 


Likewise, the comparison of the ohmic losses, representing a global quantity, shows a good match: The 2D model outputs a length-related power loss of $ 21.9 \,\mu\mathrm{Wm}^{-1}$, whereas the 3D model has a total loss of $ 4.34 \,\mu\mathrm{W}$. Scaling the length-related losses up to the cable's length results into a loss of $ 4.38 \,\mu\mathrm{W}$ which deviates $0.9\,\%$ from the 3D result. We suspect that this discrepancy is mainly due to the different excitation types.


\section{Conclusion \& outlook}
\label{sec:conclusion}

Exploiting the helicoidal symmetry of power cables significantly reduces computational costs for their numerical analysis. We posed a coordinate-transformed 2D eddy current BVP on a symmetry cell, which we derived using the mathematical tool of envelopes. The BVP itself was solved using the finite element method based on the $\mathbf{H}$-$\varphi$-formulation, in which we presented a way how to achieve application desirable (co)homology spaces. The comparison with a 3D reference model showed an accurate agreement, both in local and global results. Therefore, the presented 2D model forms the basis for future work in which, e.g., more realistic cable geometries (e.g., larger number of conductors and layers) are investigated.


\section*{Acknowledgment}

Special thanks to Marcus Christian Lehmann, who provided helpful implementation advice throughout the work.

\ifCLASSOPTIONcaptionsoff
  \newpage
\fi



\end{document}